\documentclass[12pt,leqno]{article}
\usepackage[dvips]{graphicx}
\usepackage{amsmath,amssymb,amsthm,amscd}
\usepackage[mathscr]{eucal}
\usepackage{amsfonts}
\begin{document}
\parskip=6pt
\newtheorem{prop}{Proposition}
\numberwithin{prop}{section}
\newtheorem{thm}{Theorem}
\numberwithin{thm}{section}
\newtheorem{corr}{Corollary}
\numberwithin{corr}{section}
\newtheorem{lemma}{Lemma}
\newtheorem{defn}{Definition}
\numberwithin{defn}{section}
\numberwithin{lemma}{section}
\numberwithin{equation}{section}
\newcommand{\bC}{\mathbb C}
\newcommand{\bP}{\mathbb P}
\newcommand{\bN}{\mathbb N}
\newcommand{\bA}{\mathbb A}
\newcommand{\bR}{\mathbb R}
\renewcommand\qed{ }


\def\Affine {{\mathbb A}}
\def\N {{\mathcal N}}
\def\gen {{\rm gen}}

\def\P {{\mathbb P}}     
\def\C {{\mathbb C}}     


\begin{titlepage}
\title{\bf Rationally Connected Varieties and Loop Spaces\thanks
{Research 
partially supported by NSF grant DMS0203072 and OTKA grants T046378,
T42769}}
\author{L\'aszl\'o Lempert\thanks{This research was done while the first
author, 
on leave from Purdue University, visited the Department of Analysis, 
E\"otv\"os 
University, Budapest.
He is grateful to both institutions.}\\Department of  Mathematics
\\Purdue University\\West Lafayette, IN
47907-2067, USA\and Endre Szab\'o\\R\'enyi Institute of the Hungarian 
Academy of 
Sciences\\ 1364 Budapest\\ PO Box 127}
\end{titlepage}
\date{}
\maketitle
\abstract
We consider rationally connected complex projective manifolds $M$ and
show that 
their loop spaces---infinite dimensional complex manifolds---have
properties 
similar to those of $M$.
Furthermore, we give a finite dimensional application concerning
holomorphic 
vector bundles over rationally connected complex projective manifolds.
\endabstract

\setcounter{section}{-1}
\section{Introduction}

Let $M$ be a complex manifold and $r=0,1,\ldots,\infty$.
The space $C^r(S^1,M)$ of $r$ times continuously differentiable maps $x
\colon 
S^1\to M$, the (free) $C^r$ loop space of $M$, carries a natural
complex 
manifold structure, locally biholomorphic to open subsets of Banach
$(r<\infty)$ 
resp.
Fr\'echet $(r=\infty)$ spaces, see [L2].
The same is true of ``generalized loop spaces''---or mapping 
spaces--- $C^r(V,M)$, where $V$ is a compact $C^r$ manifold, 
possibly with boundary; when 
$r=0$, $V$ can be just a compact Hausdorff space.
A very general question is how complex analytical and geometrical
properties of 
$M$ and its loop spaces are related.

Our contribution to this problem mainly concerns rational connectivity.
A complex projective manifold is a complex manifold, biholomorphic to a 
connected submanifold of some projective space $\bP^n(\bC)$.
Such a manifold $M$ is called rationally connected if it contains
rational 
curves (= holomorphic images of $\bP^1(\bC)$) through any finite
collection of 
its points.
This is equivalent to requiring  that for a nonempty open $U\subset M
\times M$ 
and any $(p,q)\in U$ there should be a rational curve through $p$ and $q
$.
For the theory of rationally connected varieties see [AK, Kl, KMM].

Projective spaces, Grassmannians, and in general complex projective
manifolds 
birational to projective spaces are rationally connected.
In a sense rationally connected manifolds are the simplest manifolds; at
the 
same time, general complex projective manifolds can be studied through 
rationally connected ones by the device of maximally rationally
connected 
fibrations
[Kl, Theorem IV.5.4].

Here is a brief description of the results presented in this paper.
For more complete formulations and for background the reader is referred
to 
Section 1.
First we prove that loop spaces $C^r(S^1,M)$ of rationally connected
complex 
projective manifolds $M$ contain plenty of rational curves, but in some
other 
mapping spaces $C^r (V,M)$ rational curves are rare.
Then we shall discuss holomorphic functions and, more generally,
holomorphic 
tensor fields.
Extending earlier results of Dineen--Mellon and the first author [DM,L2]
we show 
that on mapping spaces of rationally connected complex projective
manifolds $M$ 
holomorphic functions are locally constant, and the same is true on
submanifolds 
of $C^r(V,M)$ consisting of so called based maps.
This infinite dimensional result has the finite dimensional corollary
that 
holomorphic linear connections on vector bundles $E\to M$ are trivial.

Next we consider the ``trivial'' component of $C^r(V,M)$ consisting of 
contractible maps.
We show that the constancy of holomorphic functions on this component
already 
follows once we know $M$ is compact, connected, and all contravariant
symmetric 
holomorphic tensor fields on $M$ (of positive weight) vanish (i.e., the
only 
holomorphic section the symmetric powers of $T^* M$ admit is the zero
section).
This property is probably weaker than rational connectivity.
In fact, rational connectivity implies that holomorphic contravariant
tensor 
fields, symmetric or not, vanish, see [AK, Theorem 30]; and
conjecturally the 
converse, ``Castelnuovo's criterion'', is also true for complex
projective 
manifolds.
Finally we prove that for compact connected $M$, if on $M$ all
contravariant 
holomorphic tensor fields of positive weight vanish then the same holds
on the 
trivial component of $C^r(V,M)$.

\section{Background and results}

\noindent
{\bf 1.1.\ Mapping spaces}.
Fix $r=0,1,\ldots,\infty$ and a compact manifold $V$ of class $C^r$,
possibly 
with boundary (or just a compact Hausdorff space, when $r=0$).
We start by quickly describing the complex manifold structure on the
mapping 
space $X=C^r(V,M)$ of a finite dimensional complex manifold $M$.
For generalities on infinite dimensional complex manifolds,
see [L1, Section 2].
We need

\begin{lemma}There are an open neighborhood $D\subset M\times M$ of the
diagonal 
and a $C^\infty$ diffeomorphism $F$ between $D$ and a neighborhood of
the zero 
section in TM with the following properties.
Setting $D^w=\{z\in M\colon (z,w)\in D\}$ and $F^w=F(\cdot,w)$, for all
$w\in M$ 
we have

\noindent (a)\ $F^w$ maps $D^w$ biholomorphically on a convex subset of
$T_w M$;

\noindent (b)\ $F^w(w)\in T_w M$ is the zero vector;

\noindent (c)\ $dF^w(w)\colon T_w D^w=T_w M\to T_w M$ is the identity.

(In (c) we have identified a tangent space to the vector space $T_w M$
with the 
vector space itself.)
\end{lemma}

\begin{proof}
When $M$ is a convex open subset of some $\bC^n$ so that $TM$ is
identified with 
$\bC^n\times M$, one can take $D=M\times M$ and $F(z,w)=(z-w,w)$.
A general $M$ being locally biholomorphic to convex open subsets of $
\bC^n$, one 
obtains a covering of $M$ by open sets $W$ and $C^\infty$ maps $F_W
\colon 
W\times W\to TW$ that satisfy $(a,b,c)$ for $w\in W$ (with $D^w$
replaced by 
$W$).
If $\{\chi_W\}_W$ is a corresponding $C^\infty$ partition of unity on $M
$, one 
can take as $F(z,w)$ the restriction of
$$
\sum_W\chi_W(w) F_W(z,w)
$$
to an appropriate neighborhood of the diagonal.
\end{proof}

Given $D$ and $F$, define the complex structure on $X=C^r(V,M)$ as
follows.
A coordinate neighborhood of $y\in X$ consists of those $x\in X$ for
which 
$(x(t),y(t))\in D$ for all $t\in V$.
This neighborhood is mapped to an open subset of $C^r(y^* TM)$, the
space of 
$C^r$ sections of the induced bundle $y^* TM\to V$, by the map
$$
\varphi_y\colon x\mapsto\xi,\qquad \xi(t)=F(x(t),y(t)).
$$
It is straightforward that the local charts $\varphi_y$ are
holomorphically 
related and so define a complex manifold structure on $C^r(V,M)$; this
structure is independent of the choice of $D$ and $F$.

The above construction is slightly simpler then the one in [L2, Section
2].
Its drawback is that it does not generalize to infinite dimensional
manifolds 
$M$, that may not admit $C^\infty$ partitions of unity.
By contrast, the construction in [L2] does generalize, since it uses
partitions 
of unity on $V$ only.

A closed $A\subset V$ and $x_0\in X$ determine a subspace of ``based''
maps.
Denoting the $r$--jet of $x$ by $j^r x$, the subspace in question is
\begin{equation}
Z=C^r_{A,x_0 }(V,M)=\{x\in X\colon j^r x|A=j^r x_0|A\},
\end{equation}
a complex submanifold of $X$.
As explained in [L2, Sections 2,3], for $x\in X$ the tangent space $T_x
X$ is 
naturally isomorphic to $C^r(x^*TM)$; if $x\in Z$, under this
isomorphism $T_x 
Z\subset T_x X$ corresponds to
\begin{equation}
C_{A}^r (x^* TM)=\{\xi\in C^r (x^* TM)\colon j^r\xi|A=0\}.
\end{equation}

Up to this point $TX$, $TZ$ are real vector bundles.
However, as in finite dimensions, the local charts endow the real
tangent 
bundles of $X$ and $Z$ with the structure of a locally trivial
holomorphic 
vector bundle, and we shall always regard $TX$ and $TZ$ as such.

\setcounter{thm}{2}
\medskip\noindent
{\bf 1.2.\ Rational connectivity}.
While our principal interest is in complex manifolds, we will have to
deal with 
projective (or quasiprojective) varieties defined over fields other than
$\bC$ 
as well.
Then we shall use the language of algebraic geometry, in particular the
topology 
implied will be Zariski's.
If $M$ is a variety defined over a field $k$, we write $M(k)$ for its
points 
over $k$.
When $k$ is algebraically closed we shall ignore the difference between
$M$ and 
$M(k)$, so, for instance, a smooth projective variety $M$ over $\bC$
will be 
thought of as a complex projective manifold determined by $M(\bC)$.

\setcounter{defn}{1}
\begin{defn} \label{definition.1.2}
Let $M$ be a smooth projective variety defined over a field $k$ of 
characteristic 0.
When $k$ is algebraically closed and uncountable, $M$ is rationally
connected if 
there is a morphism $f\colon\bP^1\to M$ defined
over $k$ (i.e., a rational curve) such 
that the induced subbundle $f^* TM$ is ample:\ $f^* TM\approx\bigoplus
{\cal O}_{{\mathbb P}^1}(d_j)$, with all $d_j>0$.
In general, $M$ is rationally connected if it is such when considered
over some 
(and then over an arbitrary) uncountable algebraically closed field $K
\supset 
k$.
\end{defn} 

Over a field of positive characteristic the above property defines the
so called 
separably rationally connected varieties.
When $k=\bC$, rational connectivity is equivalent to requiring that
there be a 
nonempty open $U\subset M(\bC)\times M(\bC)$ such that for $(p,q)\in U$
there is 
a rational curve through $p$ and $q$; and also to requiring that through
any 
finite collection of points in $M(\bC)$ there be a rational curve.
For all this, see [AK, Definition--Theorem 29].

\medskip\noindent
{\bf 1.3.\ Rational connectivity of loop spaces}.
\begin{thm}Let $M$ be a rationally connected
complex projective manifold and $V$ a one real dimensional manifold.
Then the space $C^r(V,M)$ is rationally connected in the sense that for
any 
$n\in \bN$ there is a dense open $O\subset C^r(V,M^n)$ such that through
any 
$n$--tuple of maps $(x_1,\ldots,x_n)\in O$ there is a rational curve in 
$C^r(V,M)$.
\end{thm}

Taking $V=S^1$ we see that $M$ must be simply connected, a result first
proved 
by Campana [C].

The theorem would not hold for higher dimensional $V$.
First, the space $C^r(V,M)$ may be disconnected, which precludes
rational 
connectivity.
But even within components rational curves will be scarce, typically.
Let us call the component of $C^r(V,M)$ containing constant maps the
trivial 
component.

\begin{thm}If $V$ is a closed connected surface and $h\colon\bP^1
(\bC)\to 
C^r(V,\bP^1(\bC))$ holomorphic then $h$ is constant or else maps into
the 
trivial component of $C^r(V,\bP^1(\bC))$.
\end{thm}

We do not know whether in the trivial component generic $n$--tuples can
be 
connected with rational curves; but at least a nonempty open set of $n
$--tuples 
can be.
We shall not prove this, but it follows along the lines of Section 2 
(with Lemma 2.2 slightly modified), even for arbitrary $V$ and the
mapping space $C^r(V,M)$ of a rationally 
connected $M$ instead of $\bP^1(\bC)$.

\medskip\noindent
{\bf 1.4.\ Holomorphic functions}.
A simple consequence of Theorem 1.3 is that on loop spaces of
rationally 
connected complex projective manifolds holomorphic functions are
constant.
It turns out that this generalizes to spaces $Z$ of based maps, see
(1.1), even 
though a typical $Z$ will not contain compact subvarieties, let alone
rational 
curves according to [L2, Theorem 3.4].

\begin{thm}If $M$ is a rationally connected complex projective
manifold, 
$A\subset V$ is closed, and $x_0\in C^r(V,M)$, then holomorphic
functions 
$C^r_{A,x_0}(V,M)\to\bC$ are locally constant.
\end{thm}

The case when $M$ is a projective space was known earlier, see [DM,
Theorems 7, 
11] for $r=0$ and [L2, Theorem 4.2] in general.
The theorem has the following

\setcounter{corr}{5}
\begin{corr}If a holomorphic vector bundle (possibly with fibers Banach
spaces) 
$E\to M$ over a rationally connected complex projective manifold admits
a 
holomorphic (linear) connection then both $E$ and the connection are
trivial.
\end{corr}

In fact, in section 4 we shall prove a rather more general result.
However, something far more general may also be true that has nothing to
do with 
rational connectivity.
We conjecture that Corollary 1.6 is true for all simply connected
compact 
K\"ahler manifolds $M$.---As Koll\'ar noted, when $E$ has finite rank
and
$M$ is not only rationally connected but Fano, Corollary 1.6 immediately
follows from [AW, Proposition 1.2]. Holomorphic connections have been
studied first by Atiyah \cite{A}. Among other things he
proved a Lefschetz type theorem on generic hyperplanes, and 
completely classified holomorphic connections over Riemann surfaces
in terms of the fundamental group.

\medskip\noindent
{\bf 1.5.\ Holomorphic tensor fields}.
All tensor fields will be contravariant, without explicit mentioning.
Over a finite dimensional manifold $N$ these tensor fields are sections
of 
tensor powers of $T^* N$.
To avoid, in the infinite dimensional case, dealing with the ambiguous
notion of 
tensor product of Banach or Fr\'echet spaces/bundles, we simply define
a 
holomorphic tensor field on a complex manifold $N$ as a holomorphic
function
$$
g\colon T^j N=TN\oplus\ldots\oplus TN\to\bC,
$$
multilinear on each fiber.
The integer $j=0,1,\ldots$ is the weight of the tensor field; a tensor
field of 
weight 0 is just a holomorphic function on $N$.
If $g$ is symmetric on the fibers we speak of a symmetric tensor field.
Examples of symmetric holomorphic tensor fields are the zero fields (for
all 
weights) and the constant fields of weight 0.
We call these fields trivial, and we shall be interested in manifolds $M
$ on 
which all holomorphic tensor fields are trivial (this implies $M$ is
connected).
As said earlier, rationally connected complex projective manifolds are
of this 
kind.

In the next theorem $M$ can be a complex manifold locally biholomorphic
to open 
sets in Banach spaces.

\setcounter{thm}{6}
\begin{thm} Let $M$ be a complex manifold and $Y\subset C^r(V,M)$
a connected neighborhood of the space of constant maps.

(a)\ If on $M$ all symmetric holomorphic tensor 
fields are trivial then holomorphic functions on $Y$
are constant.

(b)\ If on $M$ all holomorphic tensor fields are trivial then 
the same holds on $Y$.
\end{thm}

As already said, for complex projective manifolds having only trivial 
holomorphic tensor fields is conjecturally equivalent to rational
connectivity.
According to Theorems 1.3, 1.7, both properties are inherited by loop
spaces, a 
fact we consider a mild additional evidence in favor of the conjecture.

\section{Rational connectivity of loop spaces}

In this section we shall prove Theorem 1.3.
The key is the following.

\begin{lemma}  \label{lemma.2.1}
Let $k$ be a field of characteristic zero and $M$ a smooth, 
rationally connected projective variety over $k$.
Given distinct points $p_1,\ldots,p_n\in\bP^1$ defined over $k$, there
are a 
smooth variety $W$ and a morphism $f\colon\bP^1\times W\to M$ over $k$
such that 
the map
\begin{equation}
\varphi= (f(p_\nu,\cdot))^n_{\nu=1}\colon W\to M^n
\end{equation}
is a surjective submersion on a dense open $U\subset M^n$ and its fibers
are 
irreducible.
\end{lemma}


In the proof we will have to extend the field $k$. Since after 
a field extension our varieties might become reducible, and so
cease to be varieties, we shall work in the larger category of 
schemes of finite type over a field. When dealing with a base field 
other than $k$, we will indicate it in the subscript.
Note that the notion of direct 
product of varieties (or schemes) depends on the base-field.
For a good introduction to the language of schemes, 
and a detailed explanation of the basic properties of families and
their fibers we refer to chapters II/1--II/3 of \cite{Hs}. Below
we summarize what will be needed for the proof of Lemma 2.1.

\setcounter{defn}{1} 
\begin{defn}
All our schemes are assumed to be schemes of finite type without
explicitly writing so.
\end{defn}

\begin{defn}[Dominant families] 
Let $k$ be a field of characteristic 0 and $N$ a variety over $k$.
A {\em dominant family} over $N$ is a morphism $A\to N$ 
from a $k$-scheme to $N$ with dense image.
The fiber of this family over a closed point $p\in N$ is just the
inverse image of $p\in N$, a $k$-scheme.
A {\em generically defined morphism} between two dominant families 
$\alpha:A\to N$ and $\beta:B\to N$ is an equivalence class of pairs
$(U,\phi)$, where $U\subseteq N$ is an open dense subset and 
$\phi:\alpha^{-1}(U)\to B$ is a morphism such that
$\beta\circ\phi=\alpha|{\alpha^{-1}(U)}$ (i.e. $\phi$ acts
fiber-wise). 
Two pairs $(U,\phi)$ and $(V,\psi)$ are considered equivalent, 
if the restrictions of $\phi$ and $\psi$ to the subset 
$\alpha^{-1}(U\cap V)$ are equal. 
Thus the dominant families over $N$ form a category, 
whose morphisms are the generically defined morphisms.
\end{defn}

Restricting a dominant family
to  a family over a dense open $P\subset N$
defines an equivalence of categories.
So, in the following discussion we shall assume that our $N$ is an
affine variety with coordinate 
ring $\N$. Let $K$ denote the field of rational functions on $N$, 
the quotient field of $\N$.

If $A\subset N\times \Affine_k^n$ is a closed subscheme for some $n$,
where $\Affine_k^n$ denotes the $n$-dimensional affine space over $k$,
and the restriction $A\to N$ of the projection 
$N\times \Affine_k^n\to N$ is dominant, we call the family $A\to N$
{\em affine}.
Then $A$ is simply defined via the vanishing of certain $n$-variable
polynomials whose coefficients are regular functions on $N$,
i.e. elements from $\N$.
Each dominant family is the union of finitely many affine dominant
families.
A standard method for dealing with families is 
first to study the situation for affine families, 
and then to glue them together. 

Let $B\subset N\times\Affine_k^m$ be another affine
dominant family,
and let $(U,\phi):A\to B$ represent a generically defined morphism
between them. Again, we may shrink $U$ to an affine subset.
Then $\phi$ is simply given via $m$ coordinates, each coordinate an
$n$-variable polynomial whose coefficients are regular functions on
$U$, i.e., elements of $K$.

$A$ and $B$ were defined via polynomial equations with coefficients in
$\N$. Now we think of them as polynomials with coefficients in $K$,
and define the {\em generic fibers} $A_\gen\subset\Affine_K^n$
and $B_\gen\subset\Affine_K^m$ via the same equations.
These are affine $K$-schemes. 
Moreover, the $m$ polynomials we used to define $(U,\phi)$ also have
coefficients in $K$. Hence they give us a morphism
$\phi_\gen:A_\gen\to B_\gen$, 
and we call this the {\em generic fiber of $(U,\phi)$}. 
It is easy to check 
that taking generic fibers is a functor from the
category of affine dominant families defined over $N$ into the
category of affine $K$-schemes
(in particular, it does not depend on the chosen affine embeddings). 

There is also a functor in the opposite direction obtained as follows: 
Consider an arbitrary affine $K$-scheme defined
via finitely many polynomial equations with coefficients in $K$. 
Let $D\in\N$ be a common denominator for these coefficients and
let $N_D\subset N$ denote the complement of the zero set of $D$.
Then the coefficients are regular functions on $N_D$, hence our
$K$-scheme is isomorphic to the generic fiber of the affine dominant
family over $N_D$ defined via the same equations.
A similar argument proves that each morphism of affine $K$-schemes is 
just the generic fiber of a generically defined
morphism of affine dominant families. 
Hence taking generic fiber is an equivalence of categories.

The construction of generic fibers in affine dominant families being
canonical, one can easily extend it to arbitrary dominant families, 
by simply gluing together the affine pieces. Before we state the next
theorem, we recall a definition:

\begin{defn}  \label{field-extension}
Let $X$ be any $k$-scheme, and $K\supset k$ a field
extension.
Then $X$ is defined via gluing together certain affine subsets,
the affine pieces are given via polynomial equations with
coefficients in $k$. If we consider the same equations over $K$, 
we get affine $K$-schemes.
These affine pieces  glue together the same way as the original 
affine $k$-schemes to produce a $K$-scheme $X_K$.
\begin{enumerate}
\item
$X_K$ will denote the $K$-scheme we obtain.
\item
If $p\in Y\subset X$ is a $k$-point and a $k$-subscheme, then the
above procedure gives us a $K$-point and a 
$K$-subscheme
$p_K\in Y_K\subset X_K$. 
\item
We say that $X$ is {\em geometrically irreducible} if $X_K$ is
irreducible for all field extensions $K\supset k$.
\end{enumerate}
\end{defn}

We note that an open dense subset of a scheme $X$ is geometrically
irreducible precisely when $X$ is.

\setcounter {thm}{4}
\begin{thm}[Families and their fibers]  \label{generic-fiber}
Let $k$ be a field of characteristic zero, $N$ a variety over $k$,
and $K$ the field of rational functions on $N$.
\begin{enumerate}
\item  \label{generic-fiber-exists}
Each dominant family $A\to N$ has a {\em generic fiber}
$A_\gen$, which is a $K$-scheme.
Each generically defined morphism 
between families $A\to N$ and $B\to N$ 
has a generic fiber $A_\gen\to B_\gen$.
\item  \label{generic-fiber-functor}
Taking generic fiber is a functor from the category of dominant
families of varieties over $N$ to the category of $K$-schemes. It is
an equivalence of categories.
\item  \label{generic-fiber-parts}
For a $k$-scheme $X$ the projection $X\times N\to N$ is called the 
{\em trivial family with fiber $X$}.
It's generic fiber is $X_K$.
Under the above equivalence of categories the subschemes
$Y_\gen\subset X_K$ correspond to {\em families of subschemes of $X$} 
parameterized by $N$, i.e. subschemes $Y\subset X\times N$ which 
dominate $N$;
and $K$-points  $q_\gen\in X_K$ correspond to {\em families of points 
of $X$} parameterized by $N$, i.e. subvarieties $q\subset X\times N$ 
which project birationally to $N$.
Furthermore, $q_\gen\in Y_\gen$ if and only if after
restricting to some dense Zariski open subset $U\subset N$ the second
family contains the first: $q|_U\subset Y|_U$.
For a $k$-point and a $k$-subscheme $q'\in Y'\subset X$ the constant
families $q=q'\times N$, $Y=Y'\times N$ correspond to the $K$-point
$q_\gen=q'_K\in X_K$ and $K$-scheme $Y_\gen=Y_K\subset X_K$.
\item  \label{generic-fiber-product}
If $X\to N$ and $Y\to N$ are dominant families, then the fiber product 
$X\times_NY\to N$ is again a dominant family. 
Moreover, the corresponding $K$-scheme is 
$(X\times_NY)_\gen\simeq X_\gen\times_K Y_\gen$,
where $\times_K$ stands for the direct product of $K$-schemes.
\item  \label{generic-fiber-pullback}
If $f:M\to N$ is a dominant morphism of $k$-varieties with function
fields $L\supset K$ and $X\to N$ is a dominant family, 
then the {\em pullback family} is the fiber product $X\times_NM\to M$, 
a dominant family over $M$ whose generic fiber 
is isomorphic to the $L$-scheme $(X_\gen)_L$.
Moreover, the fiber of the pullback family at any $k$-point $p\in M$
is naturally isomorphic to the fiber of the original family at $f(p)$. 
\item  \label{generic-fiber-geom-irred}
If the generic fiber of a dominant family is geometrically irreducible 
then almost all fibers are irreducible, hence connected in the Zariski
topology. 
\end{enumerate}
\end{thm}

{\em Sketch of proof.}
We already discussed the first two statements;
\ref{generic-fiber-parts}, \ref{generic-fiber-product} and
\ref{generic-fiber-pullback}  follow easily from
\ref{generic-fiber-exists}, \ref{generic-fiber-functor} and the
definitions.
For a detailed explanation of the basic properties of families and
their fibers we refer to chapters II/1--II/3 of \cite{Hs}.
We give a short proof for \ref{generic-fiber-geom-irred} here.
Supposing that $\alpha:A\to N$ is a dominant family over $k$
such that $A_\gen$ is geometrically irreducible, we shall show that
its fibers are geometrically irreducible over a dense Zariski open
subset of $N$.
Our $\alpha$ is defined over a finitely generated subfield of $k$, so
we may assume that $k$ is finitely generated, hence a subfield of
$\C$, the field of complex numbers. Then we may extend $k$ to $\C$
(reducible schemes stay reducible after field extensions), 
so from now on $k=\C$.
We may replace $A$ with its reduced scheme structure.
Since the regular part of $A$ intersects almost all fibers in a dense
Zariski open subset, we may assume that $A$ is a smooth variety.

Next we pick an irreducible subvariety $S\subset A$ of dimension
$\dim N$ such that the restricted map $\alpha|_S:S\to N$ is dominant
(if $N$ is at least one-dimensional then we may simply intersect $A$
with $\dim A-\dim N$ hyperplanes in general position, if $N$ is a
point then take any closed point for $S$).
We shall think of $S$ as a multiple valued section of $A\to N$.
By the algebraic version of Sard's lemma
(III/10.7 in \cite{Hs}) 
we can find a dense, Zariski open subset $U\subset N$ such that the
restrictions  
$\alpha^{-1}(U)\to U$ and $(\alpha|_S)^{-1}(U)\to N$ are submersions,
and further shrinking $U$ we can achieve that the second map is
finite, hence a covering space in the topological sense.
For simplicity we replace $N$ with $U$, so from now on $\alpha$ is a
submersion, and the restriction $S\to N$ is a finite covering space.
We may do this since $U$ is dense, hence we kept almost all fibers of
the original family.

Let $\alpha_S:A_S\to S$ denote the pullback of the family $\alpha$
via $S\to N$. 
The virtue of this family is that it has a section $s:S\to A_S$
defined as $s(p)=(p,p)\in A\times_NS$.
In each fiber of $\alpha_S$ we consider the connected component which
intersects this section, and let $A_S^*\subset A$ be their union.
It is easy to see, that $A_S^*$ is a connected component of $A_S$ (in
the classical topology, hence also in the Zariski topology).
On the other hand, each component of $A_S$ surjects onto $A$, hence
surjects onto $N$, and also on $S$ (since proper subvarieties of $S$
may not surject onto $N$).
Hence the irreducible components of $A_S$ are in
one-to-one correspondence with the irreducible components of its
generic fiber. But this generic fiber can be obtained from $A_\gen$
by extending the base-field $K$ to the function field of $S$
(see \ref{generic-fiber-pullback} of this theorem), 
so it is irreducible by assumption. Hence $A_S$ must also be
irreducible. 
Therefore $A_S=A_S^*$, and all fibers of $\alpha_S$ are irreducible.
But the fibers of $\alpha_S$ are the same as the fibers of $\alpha$ 
(just they appear many times), we proved our claim.

\vskip 5pt
{\em Proof of Lemma \ref{lemma.2.1}.}
At the price of replacing $M$ by $\P^3\times M$ it can be assumed that
$\dim(M)\ge 3$. Our plan is to study the family of all $n$-tuples of
points on $M$. This family is parameterized by $M^n$, and it is given
by $n$ families of points of $M$:
$$
\begin{array}{cccl}
q_\nu &\hookrightarrow &M\times M^n \cr
\gamma_\nu\Big\downarrow\ \ \ \ &&
     \Big\downarrow &\kern 20pt (\nu=1,2,\dots n)\cr
M^n &=== &M^n
\end{array}
$$
where the left hand side is
$$
q_\nu = 
\Big\{(m\ ;\ m_1,m_2,\dots, m_n)\in M\times M^n\ :\ m=m_\nu\Big\}
\simeq M^n
$$ 
and
$$
\gamma_\nu(m_\nu ;\ m_1,m_2,\dots, m_n) = (m_1,m_2,\dots, m_n).
$$
Taking generic fibers we get an $n$-tuple of $K$-points
$$
q_{1,\gen},\ q_{2,\gen},\ \dots ,q_{n,\gen}\in M_K.
$$

Our next goal is to find a smooth, geometrically irreducible
$K$-variety $V$ together with a family $g:\P^1_K\times_K V\to M_K$ 
of smooth rational curves, each passing through all of the points
$q_{1,\gen},q_{2,\gen},\dots, q_{n,\gen}$.
Since $M$ is rationally connected, so is $M_K$,
cf. Definition~\ref{definition.1.2}. 
By Theorem~16 of \cite{KS}, over $K$ there is a
family of smooth rational curves in $M_K$, parameterized by a smooth,
geometrically irreducible $K$-variety $W$, 
such that each curve passes through our $n$ points.
But this theorem gives the rational curves as subvarieties of $M_K$,
so our $V$ will be a $PGL_K(2)$-bundle over $W$ to account for all
parameterizations. 
(Note here that although we found an entire family of curves,
$V$ might not have any point defined over $K$,
so we cannot easily get a single rational curve defined over $K$.)

By the smoothness assumption each rational curve $g(\cdot,v)$ passes
through all $q_{\nu,\gen}$ exactly once, whence there are
$K$-morphisms
$$
\sigma_\nu:V\to\P^1_K\ ,\kern 15pt\nu=1,2,\dots n
$$
such that $g(\sigma_\nu(v),v)=q_{\nu,\gen}$ for all $v\in V$. Indeed,
for each $\nu$ the projection
$$
S=\Big\{(s,v)\in\P^1_K\times_K V : g(s,v)=q_{\nu,\gen}\Big\} \to V
$$
is a bijective morphism between smooth varieties. Its inverse
is rational by Sard's lemma, [Hs, III/10.7], hence regular
by Zariski's Main Theorem, [Hs, V/5.2]. Composing this inverse 
with the projection
$S\to\P^1_K$ we obtain $\sigma_\nu$. By shrinking $V$ to a dense open
subset we may achieve that $\sigma_\nu(v)$ is never
$\infty\in\P^1_K$. 

Next we can think of the given
points $p_\nu\in\P^1_k$ as points $p_{\nu,K}\in\P^1_K$ with the same
homogeneous coordinates (see Definition~\ref{field-extension}).
We can compose $g(\cdot,v)$ with appropriate interpolating
polynomials---or rational functions if some $p_\nu=\infty$---to
replace $\sigma_\nu(v)$ with $p_\nu$ and arrive at the situation
\begin{equation} \label{g-p-q}
g(p_{\nu,K},v) = q_{\nu,\gen}\ ,\kern15pt
v\in V,\ \nu=1,2,\dots n.
\end{equation}
The rational curves $g(\cdot,v)$ may no longer be smooth, but this will
not matter.

The $K$-variety $V$ is the generic fiber of a dominant family
$\phi:V^*\to M^n$ whose fibers are reduced (since $V$ is reduced).
By Theorem~\ref{generic-fiber} we may assume that each fiber of $\phi$
is irreducible (we simply shrink $V^*$ to a dense open subset). 
Upon further shrinking, by Sard's lemma (III/10.7 in \cite{Hs}) 
we may also assume that $\phi$ is a submersion. 
We use again Theorem~\ref{generic-fiber}, and see that
$\P^1_K\times_KV$ is the generic fiber of the dominant family
$(\P^1\times M^n)\times_{M^n}V^*\simeq \P^1\times V^*\to M^n$
and $g$ is the generic fiber of a generically defined morphism of
families:
$$
\begin{array}{ccccc}
\P^1\times V^*&\stackrel{g^*}\dashrightarrow &M\times M^n\cr
\kern 38pt\Bigg\downarrow\phi\circ{\rm pr}_{V^*} &&\Big\downarrow\cr
M^n &=== &M^n
\end{array}
$$
Since $g^*$ is a rational map into a projective variety, it
extends (uniquely) to the complement of a codimension 2 subvariety of
$\P^1\times V^*$. Hence we may shrink $V^*$ further, and achieve that
$g^*$ is an everywhere defined morphism 
(see the paragraph about fundamental points on page 50 of \cite{M}).
By Theorem~\ref{generic-fiber}, equation~(\ref{g-p-q}) implies that
$$
g^*(p_\nu,v) = \gamma_\nu^{-1}(\phi(v)) = \Big(\phi(v)_\nu;\ \phi
(v)\Big)
\in M\times M^n,
$$
where for $v\in V^*$, $\phi(v)_\nu\in M$ denotes the $\nu$-th coordinate of
$\phi(v)\in M^n$, $\nu=1,\dots,n$.
Now we set $W=V^*$, and $f={\rm pr}_M\circ g^*:\P^1\times V^*\to M$, 
the first component of $g^*$. Then $f(p_\nu,v)=\phi(v)_\nu$, hence
$$
\Big(f(p_\nu,\cdot)\Big)_{\nu=1}^n = \phi:W^*\to M^n
$$
is a submersion, and its fibers are irreducible.
This proves the lemma.


Next we need a result from differential geometry.
Let $V$ be a one dimensional compact manifold.

\begin{lemma}Let $\varphi\colon W\to U$ be a surjective $C^\infty$
submersion 
between finite dimensional $C^\infty$ differential manifolds, whose
fibers are 
connected.
For any $r$ and $y\in C^r (V,U)$ there is such an $\eta\in C^r(V,W)$
that 
$\varphi\circ \eta=y$.
\end{lemma}

\noindent
{\it Proof}.
First observe that any compact subset $C$ of a fiber $\varphi^{-1}(u)$
has an 
open neighborhood $W_0\subset W$ such that $\varphi|W_0$ is a trivial
fiber 
bundle with connected fibers.
To verify this we can assume $U=\bR^m$ and $u=0$.
A partition of unity argument gives a connection on $W$, i.e.~a
subbundle 
$H\subset TW$ complementary to the tangent spaces of the fibers of $
\varphi$.
Fix a relatively compact, connected open neighborhood 
$G\subset\varphi^{-1}(u)$ of $C$.
Connect an arbitrary $v\in\bR^m$ with $0\in\bR^m$ by a curve $\gamma$
consisting 
of $m$ segments, the $\mu$'th segment parallel to the $\mu$'th
coordinate axis.
If $v$ is in a sufficiently small neighborhood $U_0\subset\bR^m$ of $u$
and 
$c\in G$ then $\gamma$ can be uniquely lifted to a piecewise smooth
curve 
$\Gamma$, tangent to $H$ and starting at $c$.
Let $\psi(c,v)$ denote the endpoint of $\Gamma$.
Then $\psi$ is a fiberwise diffeomorphism of $G\times U_0$ on an open 
neighborhood $W_0$ of $C$, as claimed.

It follows that there are closed arcs $A_1,\ldots,A_n$ covering $V$ and
$C^r$ 
maps $\eta_\nu\colon A_\nu\to W$ such that $\varphi\circ\eta_\nu=y$.
We show that there is a $C^r$ map $\overline\eta\colon A_1\cup A_2\to W$
such 
that $\varphi\circ\overline\eta=y$.
Indeed, $A_1\cap A_2$ is empty or consists of one or two components.
In the first case $\overline\eta=\eta_\nu$ on $A_\nu,\ \nu=1,2$, will
do.
Otherwise choose points $b_i$ from each component of $A_1\cap A_2$;
thus 
$A_1\cup A_2\backslash \{b_i\}_i$ is the disjoint union of two arcs $
\tilde 
A_\nu\subset A_\nu,\ \nu=1,2$.
Using the neighborhoods of $C=C_i=\{\eta_1 (b_i),\ \eta_2(b_i)\}$ from
our 
initial observation, it is straightforward to construct the required 
$\overline\eta$; it will agree with $\eta_\nu$ on $\tilde A_\nu$, away
from a 
small neighborhood of $b_i$.

Now one can continue in the same spirit, fusing more and more arcs,
eventually 
to obtain the $\eta$ of the lemma.

\bigskip\noindent
{\it Proof of Theorem 1.3}.
Fix distinct $p_1,\ldots,p_n\in\bP^1(\bC)$ and apply Lemma 2.1, with $k=
\bC$.
We obtain a holomorphic map
$f\colon\bP^1(\bC)\times W\to M$ of complex manifolds so that
$$
\varphi= (f(p_\nu,\cdot))^n_{\nu=1}\colon W\to M^n
$$
is a surjective submersion on a Zariski dense open $U\subset M^n$, with 
irreducible, hence connected fibers (see [M, 4.16 Corollary]).
Since the complement of $U$ is of real codimension 2 in
$M^n$, $O=C^r(V,U)$ is dense in $C^r(V,M^n)$.
Given $x_1,\ldots,x_n\in C^r(V,M)$ such that $y=(x_1,\ldots,x_n)\in O$, 
there is a holomorphic map 
$h\colon\bP^1(\bC)\to C^r(V,M)$ with $h(p_\nu)=x_\nu,\ \nu=1,\ldots,n$.
Indeed, using Lemma 2.2 one finds $\eta\colon V\to W$ such that 
$\varphi\circ\eta=y$.
Setting
$$
F(p,t)=f(p,\eta(t)),\quad p\in\bP^1(\bC),\ t\in V,
$$
the map $h$ given by $h(p)=F(p,\cdot)$ will do.

\section{The Proof of Theorem 1.4}

Fix $d\in\bN$ and consider the space of pairs of complex polynomials
$$
\biggl\{\biggl(\sum^d_{j=0}\alpha_j u^j,\ \sum^d_{j=0}\beta_j u^j
\biggr)\colon 
(\alpha_d,\beta_d)\not= (0,0)\biggr\},
$$
a complex manifold on which the group $\bC^*$ of nonzero complex numbers
acts 
holomorphically and freely by coefficientwise multiplication.
Denote the quotient manifold by $E$, and by $\pi\colon E\to\bP^1(\bC)$
the 
projection
$$
\pi(\sum\alpha_j u^j\colon \sum\beta_j u^j)=(\alpha_d\colon\beta_d).
$$
Thus $E$ is a locally trivial fiber bundle.
Let $\Delta\subset E$ be the discriminant set, corresponding to pairs
of 
polynomials with a common zero.

\begin{lemma}If $V$ is a closed surface and $\psi\colon V\to E\backslash
\Delta$ 
is continuous then $\pi\circ\psi\colon V\to\bP^1(\bC)$ is homotopic to
a 
constant.
\end{lemma}

\noindent
{\it Proof}.
First observe that the fiber maps
$$
H_s\biggl(\sum_0^d\alpha_j u^j\colon\sum^d_0\beta_j u^j\biggr)
=\biggl(\alpha_d u^d+s\sum_0^{d-1}\alpha_j u^j\colon\beta_d 
u^d+s\sum^{d-1}_0\beta_j u^j\biggr),
$$
$1\geq s\geq 0$, deform $E$ on the image of a section.
It follows that any section is a homotopy inverse of $\pi$, in
particular the 
section
$$
\sigma(x_0\colon x_1)=(x_0 u^d-x_1\colon x_1 u^d-x_0),\qquad (x_0\colon 
x_1)\in\bP_1(\bC).
$$
Next let $L\to E$ denote the holomorphic line bundle determined by the 
hypersurface $\Delta$.
Thus $L$ has a holomorphic section that vanishes precisely on $\Delta$; 
in particular $\psi^*L$ is trivial.
One checks that the graph of $\sigma$ intersects $\Delta$ in two points,
whence 
the holomorphic line bundle $\sigma^* L\to\bP^1(\bC)$ has a holomorphic
section 
with two zeros.
It follows that $\sigma^* L$ is not (even topologically) trivial.
On the other hand, since $\sigma\circ\pi\simeq$ id$_E$,
$$
(\pi\circ\psi)^* \sigma^* L=(\sigma\circ\pi\circ\psi)^* L\approx \psi^*
L,
$$
is trivial.
Comparing Chern classes we find that $\pi\circ\psi$ must induce the zero
map 
$H^2(\bP^1(\bC))\to H^2 (V)$, and so $\pi\circ\psi$ is homotopic to a
constant 
by Hopf's theorem, see [S, Chapter 8, Section 1].

\bigskip\noindent
{\it Proof of Theorem 1.4}.
Define $g\in C^r(\bP^1(\bC)\times V, \bP^1(\bC))$ by
$$
g(s,t)=h(s)(t),\qquad s\in\bP^1(\bC),\ t\in V.
$$
For each $t\in V$ $g(\cdot,t)$ is a holomorphic map $\bP^1(\bC)\to\bP^1
(\bC)$, 
whose degree $d$ is independent of $t$; we assume $h$ is nonconstant so
that 
$d>0$.
A degree $d$ self map of $\bP^1(\bC)$ is a rational function
$$
\sum^d_{j=0}\alpha_j s^j/\sum^d_{j=0}\beta_j s^j,\qquad 
s\in\bC\subset\bP^1(\bC),
$$
with coprime numerator and denominator, and $(\alpha_d,\beta_d)\not=
(0,0)$.
Since numerator and denominator are determined up to a common factor 
$\lambda\in\bC^*$, degree $d$ maps correspond to points in $E\backslash
\Delta$, 
and $h$ induces a $C^r$ map $\psi\colon V\to E\backslash\Delta$.
By Lemma 3.1 $h(\infty)=\pi\circ\psi\in C^r(V,\bP^1(\bC))$ is
homotopically 
trivial, and therefore so are all maps $h(s),\ s\in\bP^1(\bC)$, q.e.d.

\section{Holomorphic Functions on the Manifold of Based Loops}

In this section we shall consider a rationally connected complex
projective 
manifold $M$, the space $C^r_{A,x_0}(V,M)=Z$ of based maps, $A\subset V,
\ x_0\in 
C^r(V,M)$, and we shall show that complex valued holomorphic functions
on $Z$ 
are locally constant, Theorem 1.5.
We shall also derive Corollary 1.6, in a more general form.

\begin{lemma}Given $p\in M$ and $v\in T_p M$, there are a neighborhood
$U$ of 
$p$ and a holomorphic map $\varphi\colon\bP^1(\bC)\times U\to M$ such
that
$$
\varphi(\infty,\cdot)=id_U\quad\text{and}\quad 
\varphi_* T_{(\infty,p)}(\bP^1(\bC)\times \{p\})\ni v.
$$
\end{lemma}

\noindent
{\it Proof}.
By [AK, Definition--Theorem 29] there is a  holomorphic map
$\varphi_0\colon\bP^1(\bC)\to M$ such that $\varphi_0(\infty)=p$ 
and the induced bundle $\varphi_0^* TM$ is ample. We shall obtain $
\varphi$
by deforming $\varphi_0$.
Let $\zeta$ denote the zero section of $T\bP^1(\bC)$.
By deformation theory, there are a pointed complex manifold $(W,o)$ and
a
holomorphic map $\psi\colon\bP^1(\bC)\times W\to M$ such that 
$\psi(\cdot,o)=\varphi_0$ and
\begin{equation}
T_oW\ni\omega\mapsto\psi_*(\zeta,\omega)\in{\cal O}(\varphi_0^* TM)
\end{equation}
is an isomorphism; here $(\zeta,\omega)$ is a section of 
$T(\bP^1(\bC)\times W)$ over $\bP^1(\bC)\times \{o\}$.
This follows from [Kd], that studies deformations of submanifolds of 
a given manifold.
To make the connection with deformations of maps $\bP^1(\bC)\to M$ 
needed here, observe that deformations of $\varphi_0$ correspond to 
deformations of the graph of $\varphi_0$ as a submanifold of 
$\bP_1(\bC)\times M$.

Since $\varphi_0^*TM$ is spanned by global sections, (4.1) shows
$\psi(\infty,\cdot)$ is a submersion near $o\in W$.
Shrinking $W$ we can therefore arrange that $\psi(\infty,\cdot)$ is
everywhere submersive, and so
$$
W_\infty=\{w\in W\, :\,\psi(\infty,w)=p\}
$$
is a submanifold. Under the isomorphism (4.1) $T_oW_\infty\subset T_oW$
corresponds to
$$
{\cal O}_\infty(\varphi_0^*TM)=
\{\sigma\in{\cal O}(\varphi_0^*TM)\,:\,\sigma(\infty)=0\}.
$$
We let $s$ denote the standard complex coordinate on 
$\bC\subset\bP^1(\bC)$, and for $\sigma\in {\cal O}_\infty
(\varphi_0^*TM)$
define
$$
\sigma'(\infty)=\lim_{s\to\infty}s\sigma(s)\in T_pM.
$$
As $\varphi_0^*TM$ is ample, the map $\sigma\mapsto\sigma'(\infty)$
is onto.

Note that the vector field $s^2\partial/\partial s$ extends to all
of $\bP^1(\bC)$. With $\varphi_w=\psi(\cdot,w)$ consider the map
\begin{equation}
W_\infty\ni w\mapsto\varphi_{w*}(s^2\partial/\partial s|_{s=\infty})\in
T_pM.
\end{equation}
We shall assume that (4.2) maps $o\in W_\infty$ to $0$, which
can always be arranged upon replacing $\varphi_0(s)$ by $\varphi_0
(s^2)$.
One then computes, say in local coordinates, that the differential of
the map (4.2) is the composition of (4.1) with the map
$$
{\cal O}_\infty(\varphi_0^*TM)\ni\sigma\mapsto\sigma'(\infty)\in
T_pM\approx T_0(T_pM).
$$
It follows that (4.2) is a submersion near $o$; hence
$v\in\varphi_{w*}T_\infty\bP^1(\bC)$ for some $w\in W_\infty$.
In some neighborhood $U\subset M$ of $p$ the submersion
$\psi(\infty,\cdot)$ has a holomorphic right inverse $\rho\colon U\to W
$ 
with $\rho(p)=w$; then
$$
\varphi(s,u)=\psi(s,\rho(u)),\quad s\in\bP^1(\bC),\, u\in U
$$
defines the map sought.

\bigskip\noindent
{\it Proof of Theorem 1.5}.
Let $f\colon Z\to\bC$ be holomorphic; we have to prove $df(\xi)=0$ for
all $x\in 
Z$ and $\xi\in T_x Z$.
Fix $x$.
Given $\tau\in V$ and nonzero $v\in T_{x(\tau)} M$, construct $U$ and $
\varphi$ 
as in Lemma 4.1, and with a sufficiently small neighborhood $B\subset V$
of 
$\tau$ define a $C^r$ map
$$
\Phi\colon\bP^1(\bC)\times B\ni (s,t)\mapsto\varphi(s,x(t))\in M,
$$
holomorphic in $s$.
Note that
$$
\Phi (\infty,t)=x(t)\quad ,\quad t\in B.
$$
We take $B$ compact and (when $r\geq 1$) a $C^r$ manifold with boundary.
We also arrange that $\Phi^t=\Phi(\cdot,t)$ is an immersion near $\infty
$, when 
$t\in B$.

First suppose that $\xi\in T_x Z\approx C^r_A (x^* TM)$ is supported in
the 
(relative) interior of $B$, and
\begin{equation}
\xi(t)\in\Phi^t_* T_\infty\bP^1(\bC),\quad\text{for all }t\in B.
\end{equation}
To show that $df(\xi)=0$, consider the map
$\nu\colon C^r_{(A\cap B)\cup\partial B,\infty}(B,\bP^1(\bC))\to Z$,
$$
\nu(y)(t)=
   \begin{cases} 
     \Phi(y(t),t),&\text{if $t\in B$}\\
     x(t),&\text{if $t\in V\backslash B$}
   \end{cases}
$$
By [L2, Propositions 2.3, 3.1] $\nu$ is holomorphic, and so is
$$
f\circ\nu\colon C^r_{(A\cap B)\cup\partial B,\infty}(B,\bP^1(\bC))\to
\bC.
$$
Therefore $f\circ\nu$ is locally constant by [L2, Theorem 4.2]; for the
case 
$r=0$, see the earlier [DM].
Now $\xi$ is in the range of $\nu_*$; indeed, $\xi=\nu_*\eta$, if $\eta
(t)\in 
T_\infty\bP^1(\bC)$ is defined by $\eta(t)=0$ when $t\in V\backslash B$
and 
$\Phi_*^t\eta(t)=\xi(t)$ when $t\in B$, cf.~(4.3).
It follows that
$$
df(\xi)=d(f\circ\nu)(\eta)=0.
$$

Next choose a basis $v=v_1,\ldots,v_m$ of $T_{x(\tau)}M$ and construct 
corresponding maps
$$
\Phi_1=\Phi,\Phi_2,\ldots,\Phi_m\colon\bP^1(\bC)\times B\to M.
$$
If $B$ is sufficiently small then
$$
T_{x(t)}M=\bigoplus_j \Phi^t_{j*} T_\infty \bP^1(\bC)\quad ,\quad t\in
B.
$$
For each $j$, if 
$\xi_j\in T_x Z\approx C^r_A(x^* TM)$ has support in int\,$B$ and 
satisfies (4.3), with $j$ appended, then $df(\xi_j)=0$.
Since any $\xi\in C_A^r (x^* TM)$ supported in int\,$B$ is the sum of
such 
$\xi_j$'s, we conclude each $\tau\in V$ has a neighborhood $B$ so that 
$df(\xi)=0$ when supp $\xi\subset$ int $B$.
But then a partition of unity gives $df(\xi)=0$ for all $\xi\in T_x Z$,
as 
needed.

\medskip
We shall apply Theorem 1.5 to study holomorphic connections in the
following 
setting.
Let $\pi\colon E\to N$ be a holomorphic map of complex manifolds
locally 
biholomorphic to open subsets of Banach spaces.
Assume $\pi$ is a submersion, i.e.~$\pi_*(e)\colon T_e E\to T_{\pi(e)}N$
is 
surjective for all $e\in E$.
A holomorphic connection on $E$ (or on $\pi$) is a holomorphic
subbundle 
$D\subset TE$ such that $D_e$ is complementary to Ker $\pi_*(e),\ e\in E
$.
The connection is complete if curves in $N$ can be lifted to horizontal
curves 
in $E$, i.e., for any $x\in C^1([0,1],N)$ and 
$e\in\pi^{-1}(x(0))$
there is a $y\in C^1([0,1],E)$ such that $y(0)=e$, $\pi\circ y=x$, and
$y'(t)\in 
D_{y(t)}$ for all $0\leq t\leq 1$.
The lift is unique by the uniqueness theorem for ODE's.
For example, linear connections on Banach bundles and $G$--invariant
connections 
on principal $G$ bundles---$G$ a Banach--Lie group---are complete.

The simplest example of a connection is on a trivial bundle $\pi\colon
E=F\times 
N\to N$, with $D_{(f,n)}=T_{(f,n)}(\{f\}\times N)$.
Connections isomorphic to such a connection are called trivial.
Corollary 1.6 follows from
\setcounter{thm}{1}
\begin{thm}Let $M$ be a rationally connected smooth complex projective
manifold, 
$E$ a complex manifold locally biholomorphic to open subsets of Banach
spaces, 
and $\pi\colon E\to M$ a holomorphic submersion such that on each fiber 
holomorphic functions separate points.
If $\pi$ admits a complete holomorphic connection $D$ then the
connection is 
trivial.
\end{thm}

\noindent
{\it Proof}.
The mapping space $C^1([0,1],E)$ has a natural structure of a complex 
manifold---the construction in [L2, Section 2] carries over to Banach 
manifolds.
Horizontal lift defines a map $\Lambda$ of the manifold 
$$
\{(e,x)\in E\times C^1([0,1], M)\colon\pi (e)=x(0)\}
$$
into $C^1([0,1],E)$.
This map is holomorphic.
To see this, note that for $(e,x)$ in a small neighborhood of a fixed 
$(e_0,x_0)$, and for small $\tau\in (0,1]$, finding $y=\Lambda(e,x)$
over the 
interval $[0,\tau]$ amounts to solving an ODE.
Doing this by the standard iterative scheme of Picard--Lindel\"of (see
[Hm, 
p.~8]) shows the local lift $y|[0,\tau]\in C^1 ([0,\tau],E)$ depends 
holomorphically on $(e,x)$.
Since the full lift $y$ is obtained by concatenating local lifts, $
\Lambda$ is 
indeed holomorphic.
It is also equivariant with respect to reparametrizations:\ if $\sigma
\colon 
[0,1]\to [0,1]$ is a $C^1$ map, $\sigma(0)=0$, then
\begin{equation}
\Lambda(e,x)\circ\sigma=\Lambda(e,x\circ\sigma),\qquad x\in C^1([0,1],
M).
\end{equation}

With fixed $p\in M$ and variable $q\in M$ consider
$$
Y=\{x\in C^1([0,1], M)\colon x(0)=p\},\ Y_q=\{x\in Y\colon x(1)=q\},
\text{ and }
$$
$$
Z_q=\{x\in Y_q\colon x'(0)\in T_p M\text{ and }x'(1)\in T_q M \text{ are
both 
zero}\},
$$
connected manifolds since $M$ is simply connected by [C],
or by our Subsection 1.3.
Therefore Theorem 1.5 implies that $\bC$--valued holomorphic functions
on $Z_q$ 
are constant.
In particular, for any $e\in\pi^{-1}(p)$ and holomorphic function 
$h\colon\pi^{-1}(q)\to\bC$, $h(\Lambda(e,x)(1))$ is independent of $x\in
Z_q$.
Since holomorphic functions separate points of $\pi^{-1}(q)$,
$\Lambda(e,x)(1)$ itself is independent of $x\in Z_q$.
It follows from (4.4) that $\Lambda(e,x)(1)$ is even independent of $x
\in Y_q$ 
(take e.g. $\sigma(t)=3t^2-2t^3$,
then $x\circ\sigma\in Z_q$), and so there is a holomorphic map 
$\Psi\colon\pi^{-1}(p)\times M\to E$ such that
\begin{equation}
\Lambda (e,x)(1)=\Psi (e,x(1)).
\end{equation}
One checks that $\Psi$ is biholomorphic and maps $\pi^{-1}(p)\times \{q
\}$ to 
$\pi^{-1} (q)$, $q\in M$.

To conclude, note that with $\tau\in [0,1]$ and $\sigma(t)=\tau t$
(4.4),
(4.5) imply
$$
\Lambda(e,x)(\tau)=\Psi(e,x(\tau)),
$$
i.e.~$\Psi$ maps curves $(e,x)$ to horizontal curves in $E$.
It follows that the induced connection $\Psi_*^{-1}D$ on the bundle 
$\pi^{-1}(p)\times M\to M$ is trivial, hence so is $D$.

\section{Holomorphic Tensor Fields}
To prove Theorem 1.7 we first discuss the notion of order of vanishing.
Let $Y$ be a complex manifold, locally biholomorphic to open sets in
Banach or 
even Fr\'echet spaces, $y\in Y$, and $f\colon Y\to\bC$ holomorphic.
We say that $f$ vanishes at $y$ to order $n$ if for arbitrary $0\leq k
<n$ and 
vector fields $v_1,\ldots,v_k$ on $Y$, holomorphic near $y$
$$
(v_1 v_2\ldots v_k f)(y)=0.
$$
If $Y$ is connected and $f$ vanishes at $y$ to all orders then $f\equiv
0$.

Suppose $f$ vanishes at $y$ to order $n$.
To see if it vanishes to order $n+1$, one is led to consider holomorphic
vector 
fields $v_1,\ldots,v_n$ in a neighborhood of $y$ and
\begin{equation}
(v_1 v_2\ldots v_n f)(y).
\end{equation}
Observe first that (5.1) is independent of the order in which the vector
fields 
are applied (since e.g.
$$
v_2 v_1 v_3\ldots v_n f=v_1 v_2\ldots v_n f-[v_1,v_2] v_3\ldots v_n f\\
=v_1 v_2\ldots v_n f
$$
at $y$); next that (5.1) vanishes if some $v_i$ vanishes at $y$ (since
this is 
clearly so if $v_1(y)=0$).
It follows that (5.1) depends only on the values that the $v_i$ take at
$y$, and 
so (5.1) induces a symmetric $n$--linear map
$$
d^n f(y)\colon T_y^n Y=T_y Y\oplus\ldots\oplus T_y T\to\bC.
$$

\bigskip\noindent
{\it Proof of Theorem 1.7(a)}.
Constant maps $V\to M$ form a submanifold of $Y$, biholomorphic to $M$;
we shall 
simply denote this manifold by $M\subset Y$.
If $f\colon Y\to\bC$ is holomorphic then by assumption $f|M$ is
constant.
At the price of subtracting this constant from $f$ we can assume $f$
vanishes at 
each point of $M$ to first order.
We shall prove by induction it vanishes at each $p\in M$ to arbitrary
order.

Suppose $f$ is already known to vanish to order $n\geq 1$ at each $p\in
M$, so 
that the differentials $d^n f(p)$ are defined on $T_p^n Y$.
We want to show $d^n f(p)=0$, i.e.,
\begin{equation}
d^n f(p) (\eta_1,\ldots,\eta_n)=0,\qquad \eta_i\in T_p Y,\ p\in M.
\end{equation}
Note that by Subsection 1.1 $T_p Y$ is naturally isomorphic to $C^r
(V,T_p M)$.
With fixed $\varphi_1,\ldots,\varphi_n\in C^r(V,\bC)$ define a
homomorphism 
$\Phi_n\colon T^n M\to T^n Y|M$ of holomorphic vector bundles 
\begin{equation}
\Phi_n\colon T^n M\ni (\xi_1,\ldots,\xi_n)\mapsto (\varphi_1 
\xi_1,\ldots,\varphi_n\xi_n)\in T^n Y|M;
\end{equation}
the pullback of $d^n f$ by $\Phi_n$ is a symmetric holomorphic tensor
field on 
$M$, of weight $n\geq 1$, hence vanishes.
Therefore (5.2) holds when each $\eta_i$ is of form $\varphi_i\xi_i$,
and also 
when each $\eta_i$ is a linear combination of such tangent vectors.
When $\dim T_p M<\infty$, linear combinations
$$
\sum^k_{j=1}\varphi^{(j)}\xi^{(j)},\qquad k\in\bN,\ 
\varphi^{(j)}\in C^r (V,\bC),\ \xi^{(j)}\in T_p M,\ 
$$
constitute all of $C^r(V,T_p M)$, and in general a dense subspace;
whence indeed 
$d^nf(p)=0,\ p\in M$.
This means $f$ vanishes to order $n+1$ along $M$, hence to all orders,
and 
therefore $f=0$ on $Y$. 

\medskip
For the rest of Theorem 1.7 we first extend the notions of vanishing
order and 
higher differentials to tensor fields.
Let now $f$ be a holomorphic tensor field of weight $j$ on the manifold
$Y$.
We say that $f$ vanishes at $y\in Y$ to order $n\geq 0$ if for all $0
\leq k <n$ 
and holomorphic vector fields $v_1,\ldots,v_k,w_1,\ldots,w_j$, defined
near $y$
$$
v_1\ldots v_k f(w_1,\ldots,w_j)=0\qquad\text{at }y.
$$
Note that vanishing to order 0 is automatic.
Suppose $f$ does vanish to order $n$.
As before, 
\begin{equation}
(v_1\ldots v_n f(w_1,\ldots,w_j))(y)
\end{equation}
is symmetric in the $v_l$, and for fixed $w_i$, depends only on the
values 
$v_l(y)$.

\begin{prop}If some $w_i$ vanishes at $y$ then (5.4) vanishes.
\end{prop}

\noindent
{\it Proof}.
First observe that if $F$ is a Fr\'echet space, $h$ an $F$ valued
holomorphic 
function defined in a neighborhood of $0$ in some $\bC^q$, and $h(0)=0$,
then 
there are holomorphic functions $h_1,\ldots,h_q$ such that
$$
h(z_1,\ldots,z_q)=\sum^q_{s=1} z_s h_s (z_1,\ldots,z_q)
$$
in a neighborhood of $0$.
Indeed,
$$
h(z)=\int_0^1 \frac{d}{d\lambda} h(\lambda z) d\lambda=\sum_s z_s
\int_0^1 \frac{\partial h}{\partial z_s} (\lambda z)d\lambda.
$$

Now suppose, for concreteness, that $w_1(y)=0$.
Since as far as the $v_l$ are concerned, (5.4) depends only on $v_l(y)$,
we can 
assume that all $v_l$ are tangent to a finite, say $q$, dimensional
submanifold 
$Q\subset Y$ passing through $y$.
After a local trivialization of $TY$ the above observation gives
holomorphic 
functions (local coordinates) $\zeta_1,\ldots,\zeta_q$ on $Q$ and
holomorphic 
sections $h_1,\ldots,h_q$ of $TY|Q$ near $y$, such that
$$
\zeta_1(y)=\ldots=\zeta_q (y)=0\qquad\text{and}
\qquad w_1|Q=\sum_s \zeta_s h_s.
$$
Then Leibniz's rule implies
$$
v_1\ldots v_n f(w_1,\ldots,w_j)=\sum_s v_1\ldots v_n \{\zeta_s 
f(h_s,w_2,\ldots,w_n)\}=0
$$
at $y$. 

\medskip
It follows that (5.4) depends only on ($f$ and) the values that $v_1,
\ldots,w_j$ 
take at $y$; therefore (5.4) induces a multilinear map
$$
d^n f(y)\colon T_y^{n+j} Y\to\bC.
$$

\medskip\noindent
{\it Proof of Theorem 1.7(b)}.
Assume now $f$ is a holomorphic tensor field of weight $j\geq 1$ on $Y$.
Suppose we know $f$ vanishes at all $p\in M$ to order $n\geq 0$.
As before, a choice of $\varphi_1,\ldots,\varphi_{n+j}\in C^r(V,\bC)$
defines a 
homomorphism $\Phi_{n+j}\colon T^{n+j} M\to T^{n+j} Y|M$, cf.~(5.3).
The pullback of $d^n f$ by $\Phi_{n+j}$ is a holomorphic tensor field on
$M$, 
hence $0$; from which it follows, as earlier, that $d^n f(p)=0$, $p\in M
$.
Thus $f$ vanishes to order $n+1$ along $M$, so to all orders.
This implies $f=0$ on $Y$ as claimed.

\end{document}